\documentclass[11pt]{article}

\usepackage{amsmath}
\usepackage{amsthm}
\usepackage{amsfonts}
\usepackage{amssymb}
\usepackage{amsmath,amscd}

\newcommand{\Q}{\mathbb{Q}}
\newcommand{\R}{\mathbb{R}}

\newcommand{\N}{\mathbb{N}}
\newcommand{\PP}{\mathbb{P}}
\newcommand{\RR}{\mathbb{R^{+}}}

\newcommand{\Pic}{\operatorname{Pic}}

\newcommand{\cont}{\operatorname{cont}}
\newcommand{\Exc}{\operatorname{Exc}}

\newcommand{\OO}{{\cal O}}

\newtheorem{theo}{Theorem}
\newtheorem{prop}{Proposition}
\newtheorem{lem}{Lemma}

\def\finpreuve
{\hskip 3pt \vrule height6pt width6pt depth 0pt}

\title{Classification of Fano manifolds  
containing a negative divisor isomorphic to projective space}
\author{Toru Tsukioka}
\begin{document}

\maketitle

\begin{abstract}
We classify $n$-dimensional complex Fano manifolds $X$ ($n\geq 3$)
containing a divisor $E$ isomorphic to $\PP^{n-1}$ such that $\deg N_{E/X}$ 
is strictly negative.   
\end{abstract}

\section{Introduction}
A projective manifold $X$ is called Fano manifold if its anti-canonical 
bundle $-K_{X}$ is ample.
In \cite{BCW} the authors classified $n$-dimensional 
complex Fano manifolds $X$ ($n\geq 3$) containing a divisor $E$ isomorphic to 
$\PP^{n-1}$ with normal bundle $N_{E/X}\simeq \OO_{\PP^{n-1}}(-1)$.
The purpose of the present note is to generalize their result to the case 
$\deg N_{E/X}<0$. Note that such $X$ has automatically Picard number 
$\rho(X)\geq 2$.

\begin{theo}\label{divneg}
Let $X$ be a complex Fano manifold of dimension 
$n\geq 3$ and let $E$ be a divisor of $X$.
We suppose that $E$ is isomorphic to $\PP^{n-1}$
and let $E_{|E}\simeq \OO_{\PP^{n-1}}(-d)$.
If $d>0$, the pair $(X,E)$ is (up to isomorphism) one of the following:
\footnote{For a vector space $V$, we denote by $\PP(V)$ the projective space 
of {\it lines} in $V$.}
\begin{enumerate}
\item $X\simeq\PP(\OO_{\PP^{n-1}}\oplus\OO_{\PP^{n-1}}(-d))$ ($0<d<n$)
and $E$ is the negative section such that $E_{|E}\simeq \OO_{\PP^{n-1}}(-d)$;
\item
\begin{enumerate}
\item $X$ is the blow-up of $\PP^{n}$ along $W$ 
smooth complete intersection of a hyperplane and a hypersurface of 
degree $r$ with $2\leq r \leq n$, and $E$ is 
the strict transform of the hyperplane containing $W$;
\item $X$ is the blow-up of
$\PP(\OO_{\PP^{n-1}}\oplus \OO_{\PP^{n-1}}(d'))$ 
($d'$ is an integer satisfying $-n<d'<n$)
along $W$ and $E$ 
is the strict transform of $E'$: here $E'$ 
is a section of the $\PP^{1}$-bundle 
$\PP(\OO_{\PP^{n-1}}\oplus \OO_{\PP^{n-1}}(d'))$ 
with normal bundle $\OO_{\PP^{n-1}}(d')$ and 
$W$ is a smooth divisor of degree $r$ ($d'<r<n+d'$) 
in $E'\simeq\PP^{n-1}$ 
.
\end{enumerate}
\end{enumerate}
\end{theo}

As in \cite{BCW} the proof of this theorem is based on the theory
of extremal contractions, which was originally used in \cite{MM} 
to classify Fano 3-folds with Picard number $\geq 2$. 
Recall that a fiber connected proper holomorphic map $\varphi:X\to Z$
to a normal projective variety 
 is called {\it extremal contraction} 
if $-K_{X}$ is $\varphi$-ample. We say $\varphi$ is 
{\it elemental} if $\rho(X)-\rho(Z)=1$.
By the Mori theory, a Fano manifold has a finite number of 
elementary extremal contractions. 
%and the study of these contractions 
%allows us to determine the structure of $X$ 
%(as in \cite{MM} which classifies the three dimensional smooth 
%Fano manifolds of Picard number greater than or equal to 2). 
So in our situation,  we can take 
an extremal contraction $\varphi: X \to Z$ whose 
fiber meets the negative divisor $E$. The essential 
part of the classification is to show that the restriction map  
$\varphi_{|E}:E\to \varphi(E)$ is an isomorphism. 
The corresponding fact is shown in \cite{BCW}, but their 
proof relies heavily on the fact  $\deg N_{E/X}=-1$.  
In this note we will propose another approach which allows us to get our 
result.

The argument used in this note can be applied also to classify 
other types of higher dimensional Fano manifolds of Picard number 
greater than or equal to 2. For example, the author classified 
in \cite{T}, del Pezzo surface fibrations obtained by blow-up 
along a smooth curve, in any dimension (these are automatically Fano 
manifolds of Picard number 2).

\section{Preliminary results}
Let $X$ be a complex Fano manifold of dimension $n\geq 3$ 
containing a divisor 
$E\simeq\PP^{n-1}$. We assume $\deg N_{E/X}<0$ and let 
$E_{|E}\simeq\OO_{\PP^{n-1}}(-d)$ with $d>0$. 
Note that $d<n$. Indeed, since $-K_{X}$ is ample (so is $-K_{X|E}$),  
by the adjunction formula: $K_{E}=(K_{X}+E)_{|E}$, we have
$$
0<\deg (-K_{X|E})=\deg(-K_{E})+\deg(E_{|E})=n-d.
$$ 
%Note that 
%$d\leq n-1$ by the adjunction formula: $K_{E}=(K_{X}+E)_{|E}$. 

A similar argument as 
\cite{BCW} (Lemme 1) shows that there exists an extremal ray 
$\RR[f]$ (where $f$ is a minimal rational curve of the ray) 
such that $E\cdot f>0$ and the fibers of the corresponding 
(elementary) extremal contraction $\varphi: X\to Z$ are at most of dimension 1.
By \cite{Ando}, either:
\begin{enumerate}
\item $\varphi$ is a conic bundle with smooth base $Z$, or 
\item $\varphi$ is a smooth blow-up whose centre is smooth and 
of codimension 2. 
\end{enumerate}
We shall show that the restriction map $\varphi_{|E}:E\to \varphi(E)$ 
is an isomorphism in each case.  

\begin{lem}
In the case (1), $\varphi_{|E}$ is an isomorphism. 
\end{lem}
{\bf Proof.} Since $\deg N_{E/X}<0$, by Grauert's criterion, 
there is a fiber connected holomorphic map 
$\pi : X\to Y$ to an analytic variety such that $\pi(E)$ is a point.
We have
$$
K_{X}\equiv \pi^{*}K_{Y}+\alpha E
$$
where $\alpha:=(n-d)/d$.  Here 
$\pi^{*}K_{Y}$ is to be considered as a $\Q$-divisor 
$(\pi^{*}(dK_{Y}))/d$ on $X$ 
(note that $dK_{Y}$ is Cartier divisor because $\pi(E)$ is a cyclic 
quotient singular point of order $d$). This formal expression is 
useful in the following numerical calculations.

Since the smooth projective variety $Z$ is dominated by 
$E\simeq \PP^{n-1}$ 
we have $Z\simeq \PP^{n-1}$ (see \cite{Laz}).
Hence, if we define $L:=\varphi^{*}\OO_{Z}(1)$, 
$L^{n-1}$ is numerically equivalent to 
{\it only one fiber}
of $\varphi$. Therefore,
\begin{equation*} 
\begin{cases}
L^{n}=0\\
(-K_{X})\cdot L^{n-1}=2.
\end{cases}
\end{equation*}

Since the contraction map $\varphi$ is supposed to be elemental, 
we have $\rho(X)=\rho(Z)+1=2$. Hence there exist $x,y\in\Q$ such that
$$
L\equiv x\pi^{*}(-K_{Y})-yE.
$$

We have $yd\in \N$, because $0<L\cdot e=-y(E\cdot e)=yd$ 
where $e$ is a line in $E\simeq \PP^{n-1}$.
Remark that $E^{n}=(-d)^{n-1}$ and $\pi^{*}(-K_{Y})\cdot E\equiv 0$. 
We get
$$
L^{n}=(x\pi^{*}(-K_{Y})-yE)^{n}=x^{n}m-y^{n}d^{n-1},
$$
$$
(-K_{X})\cdot L^{n-1}=(\pi^{*}(-K_{Y})-\alpha E)(x\pi^{*}(-K_{Y})-yE)^{n-1}
=x^{n-1}m-\alpha y^{n-1}d^{n-1}
$$
where $m:=(\pi^{*}(-K_{Y}))^{n}$.
Remark that $m \in \Q$ and $md^{n} \in \N$.
We have
\begin{equation*}
\begin{cases}
x^{n}m=y^{n}d^{n-1} \\
x^{n-1}m=2+\alpha y^{n-1}d^{n-1}
\end{cases}
\end{equation*}

We divide the first equality by the second:
\begin{equation}\label{x}
x=\frac{y^{n}d^{n-1}}{2+\alpha y^{n-1}d^{n-1}}.
\end{equation}
By the first equation again,
\begin{equation}\label{m}
\left(\frac{y}{x}\right)^{n}=\frac{m}{d^{n-1}}.
\end{equation}

Since
$$\frac{y}{x}=y\cdot\frac{2+\alpha y^{n-1}d^{n-1}}{y^{n}d^{n-1}},$$
we have
$$
  \left( \frac{2+\alpha y^{n-1}d^{n-1}}{y^{n-1}d^{n-1}}\right)^{n}
=\frac{m}{d^{n-1}}.
$$

We get 

$$
\left( \frac{2d^{2}}{(yd)^{n-1}}+\alpha d^{2} \right)^{n}=
\left( \frac{2+\alpha y^{n-1}d^{n-1}}{y^{n-1}d^{n-1}}\right)^{n}\cdot d^{2n}
=md^{n}\cdot d\in \N,
$$

because $md^{n}$ and $d$ are integers.
It follows that 
$$
\frac{2d^{2}}{(yd)^{n-1}}+(n-d)d\in \N 
$$ 
Hence
$l:=2d^{2}/(yd)^{n-1}$ is an integer.
We have 
$$2d^{2}=(yd)^{n-1}l. $$
On the other hand, $(l+(n-d)d)^{n}=d(md^{n})\in d\N$.
Hence $l^{n}\in d\N$.

Now we can show the equality $yd=1$: Recall first
$yd\in \N$ and $d<n$.
If $yd\geq 3$, $2(n-1)^{2}\geq 2d^{2}=(yd)^{n-1}l\geq 3^{n-1}l.$
This is impossible because $n\geq 3$.
We suppose now $yd=2$. We have $2d^{2}=2^{n-1}l$, hence
$d^{2}=2^{n-2}l$. Since $d\leq n-1$ and $l^{n}\in d\N$, this
equality is possible only when
$(n,d,l)=(3,2,2)$ or $(5,4,2)$.
These two cases can be ruled out by using the fact that
 $(-K_{X})^{n}$ is a natural number (because $X$ is smooth). Note 
first that
$$
(-K_{X})^{n}=\left(\pi^{*}(-K_{Y})-\frac{n-d}{d}E \right)^{n}
=m-\frac{(n-d)^{n}}{d}.
$$

If $(n,d,l)=(3,2,2)$ then $y=yd/d=1$, $x=1$ and $m=4$
by (\ref{x}) and (\ref{m}).
Therefore, $(-K_{X})^{3}=7/2\notin \N$, contradiction.
Similarly, if $(n,d,l)=(5,4,2)$ we have $y=yd/d=1/2$,
$x=4/3$
and $m=243/128$, so that
$(-K_{X})^{5}=211/128\notin \N$, contradiction.
Finally, we conclude that $yd=1$.

\

Now we can determine the intersection number $E\cdot L^{n-1}$ :
$$
E\cdot L^{n-1}=E\cdot(x\pi^{*}(-K_{Y})-yE)^{n-1}=(yd)^{n-1}=1.
$$
Since $L^{n-1}$ is a fiber of $\varphi$, 
this implies that the restriction map $\varphi_{|E}:E\to Z$ is an isomorphism. 
\finpreuve

\begin{lem}
In the case (2) also, $\varphi_{|E}$ is an isomorphism.
\end{lem}
{\bf Proof.} In this case, the contraction 
$\varphi : X\to Z$ is a blow-up of 
smooth center $W$ of codimension 2 and $Z$ is smooth.   
Let $F$ be the exceptional divisor of the blow-up $\varphi$.
Since $\varphi_{|F}:F \to W$ is a $\PP^{1}$-bundle,
$\tilde{W}:=(E\cap F)_{red}$ is a section (remark that
$\tilde{W}\subset F$ is contracted by 
$\pi_{|F}:F \to \pi(F)\subset Y$).

\

{\sc Claim.} {\it The intersection $E\cap F$ is transversal.}

{\it Proof:}
If not, we have an isomorphism of tangent bundles:
$TE_{|\tilde{W}}\simeq TF_{|\tilde{W}}$ 
and so $N_{\tilde{W}/E}\simeq N_{\tilde{W}/F}$.
But this is a contradiction. Indeed: $N_{\tilde{W}/E}$ is ample
(because $\tilde{W}\subset E\simeq\PP^{n-1}$) and 
 $N^{*}_{\tilde{W}/F}$ is also ample 
(because $\tilde{W}$ is contracted by $\pi_{|F}:F\to \pi(F)\subset Y$). 

\

Hence $\tilde{W}$ is a section without multiplicity, namely
$\varphi_{|E}:E \to E':=\varphi(E)$ is an isomorphism. 
In particular $E'\simeq \PP^{n-1}$. \finpreuve

\

\begin{lem}
In the case (2), $Z$ is a Fano manifold.
\end{lem}
{\bf Proof.} In fact, by the Proposition \ref{Yfano} (see below), 
it is sufficient to show that for every curve $B\subset W$, we have 
$-K_{Z}\cdot B>0$. Let $e':=\varphi_{*}e$. Since $\varphi_{|E}$ 
is an isomorphism, $e'$ is a line in $E'\simeq\PP^{n-1}$. 
Since $d<n$, we get
$$
-K_{Z}\cdot e'=-K_{X}\cdot e+F\cdot e=n-d+r>0
$$ 
where $r$ is the degree of $W$ in $E'\simeq \PP^{n-1}$.
For each curve $B$ contained in $W$ (so in $E'$), 
there exists $b\in\N$ such that $B\equiv be'$. Hence 
$-K_{Z}\cdot B=b(-K_{Z}\cdot e')>0$. \finpreuve

\

\begin{prop}\label{Yfano}
Let $X$ be a Fano manifold of dimension $n\geq 3$ and  
$\varphi:X\to Z$ a blow-up 
of center $W$ smooth subvariety of codimension 
$k$ ($2\leq k\leq n$) in a projective manifold $Z$. 
If $-K_{Z}\cdot B>0$ for any curve 
$B$ contained in $W$, then $Z$ is Fano. 
\end{prop}
{\bf Proof.} Let $F$ be the exceptional divisor of $\varphi$ 
and let $R$ be the extremal ray defining $\varphi$. 
Let $C$ be a curve of $X$ such that $[C]\not\in R$. 
If $C \nsubseteq F$, $F\cdot C\geq 0$ so that
$$
\varphi^{*}(-K_{Z})\cdot C=-K_{X}\cdot C+(k-1)F\cdot C>0.   
$$
If $C \subset F$, $\varphi_{*}C$ is an effective 1-cycle 
(because $[C]\not\in R$) 
whose support is contained in $W$. 
So, by the hypothesis of the proposition, we get  
$-K_{Z}\cdot \varphi_{*}C>0$. By \cite{Wis} Lemma (3.1), this means that 
$-K_{Z}$ is ample. \finpreuve

\section{Proof of Theorem \ref{divneg}}

In this section, we prove two propositions 
which imply our main Theorem. 

\begin{prop}
In the case (1) ($\varphi$ is a conic bundle),  we have the example 
1 in the list of Theorem \ref{divneg}.
\end{prop}
{\bf Proof.}
Since $\varphi_{|E}:E\to Z$ is an isomorphism, $\varphi$ is necessarily a 
$\PP^{1}$-bundle and
$E\cdot f=1$ where $f\simeq\PP^{1}$ 
is any fiber. 
From the exact sequence
$$
0\to\OO_{X}\to\OO_{X}(E)\to\OO_{E}(E)\to0,
$$
we deduce the exact sequence over $Z\simeq \PP^{n-1}$:
$$
0\to\OO_{\PP^{n-1}}\to\varphi_{*}\OO_{X}(E)\to\OO_{\PP^{n-1}}(-d)\to0,
$$
so that $\varphi_{*}\OO_{X}(E)$ is identified to
$\OO_{\PP^{n-1}}\oplus \OO_{\PP^{n-1}}(-d)$ and $X$ to
$\PP(\OO_{\PP^{n-1}}\oplus \OO_{\PP^{n-1}}(-d))$.
Hence we have the example 1.
\finpreuve

\begin{prop}
In the case (2) ($\varphi$ is a blow-up along a centre of codimension 2), 
we have the examples 2-(a) or 2-(b). 
\end{prop}
{\bf Proof.} 
Since $Z$ is Fano, there exists an elementary extremal contraction 
$\mu:=\cont_{\R^{+}[m]}:Z \to Z'$ such that $E'\cdot m>0$ 
where $m$ is a minimal rational curve of the ray.

\

{\it The case where there exists a fiber $M=\mu^{-1}(z')$
of dimension $\geq 2$}. In this case
there exists a curve $B\subset E'\cap M$. So $[B]\in \R^{+}[e']$
and $[B]\in \R^{+}[m]$ ($e'$ is a line in $E'\simeq \PP^{n-1}$).
Hence $\R^{+}[e']=\R^{+}[m]$. It follows that (the numerical class of) 
any curve in $E'$ is on the ray $\R^{+}[m]$, namely, $\mu(E')$ 
is a point in $Z'$. On the other hand, $E'\cdot e'>0$ because 
$E'\cdot m>0$. 
By Proposition \ref{rho}
(see below) we have $\rho(Z)=1$
and the effective divisor $E'\simeq \PP^{n-1}$ is then ample.
Finally, by \cite{BCW} Lemme 4, we conclude that 
$(Z,E')\simeq (\PP^{n},\OO_{\PP^{n}}(1))$.
In particular, $E'\cdot e'=1$.

Now we estimate the number $r:=$ the degree of $W$ as a divisor in 
$E'\simeq \PP^{n-1}$.
Since $\varphi^{*}E'=E+F$, we have $1=E'\cdot e'=E\cdot e+F\cdot e=-d+r$. 
But since
$0<d<n$, we obtain $2\leq r \leq n$.
It follows that $W=E'\cap L$ where $E'\in|\OO_{\PP^{n}}(1)|$ and
$L\in |\OO_{\PP^{n}}(r)|$ ($2\leq r \leq n$). So we get the example (2)-(a).

\

{\it The case where every fiber of $\mu$ is at most of 
dimension 1} (the following argument is essentially due to \cite{BCW}). 
By \cite{Ando}, the elementary extremal 
contraction $\mu : Z \to Z'$ is one of the following:
\begin{enumerate}
\item a $\PP^{1}$-bundle,
\item a conic bundle (with singular fibers),
\item a smooth blow-up whose centre is smooth and of codimension 2.
\end{enumerate}
Note first that only the first case is possible.
In fact in the other cases, $K_{Z}\cdot m=-1$ and there exists 
a minimal rational curve $m$ such that $m\cap W$
is not empty. If $\tilde{m}$ is the strict transform of $m$, we have
$$
K_{X}\cdot \tilde{m}=K_{Z}\cdot m+F\cdot \tilde{m}=-1+F\cdot \tilde{m}\geq 0,
$$ 
which is a contradiction because $X$ is Fano.

So, it is sufficient to study the case where $\mu$ is a
$\PP^{1}$-bundle. In this case each fiber ($=m$) meets $W$
transversally at a single point: in fact, if not, 
there exists $m$ such that $F\cdot \tilde{m}\geq 2$.
Then, as above,
$K_{X}\cdot \tilde{m}=K_{Z}\cdot m+F\cdot \tilde{m}=-2+F\cdot \tilde{m}\geq 0$,
contradiction. So $\mu_{|W}:W \to \mu(W)\subset Z'$ is an isomorphism
(by \cite{Laz}, $Z'$ is isomorphic to $\PP^{n-1}$ because $Z'$ is dominated by 
$E'\simeq\PP^{n-1}$).
We have a standard commutative diagram of Picard groups:
$$
\begin{CD}
\Pic(Z')@>>> \Pic(E')\\
@VVV         @VVV\\
\Pic(\mu(W))@>\simeq>> \Pic(W)
\end{CD}
$$
By a theorem of Lefchetz, the two vertical maps 
(induced by inclusions: $\mu(W)\subset Z'$ and $W\subset E'$) 
are isomorphisms.
Therefore $\mu^{*} : \Pic(Z')\to \Pic(E')$ is an isomorphism, so is
$\mu_{|E'} : E' \to Z'$.

Let $E'_{|E'}\simeq \OO_{\PP^{n-1}}(d')$ ($d'=-d+r$ where 
$r$ is the degree of $W\subset 
E'\simeq \PP^{n-1}$, so $d'$ might be positive).
Finally $Z\simeq \PP(\OO_{\PP^{n-1}}\oplus \OO_{\PP^{n-1}}(d'))$
and $E'$ is a section with normal bundle 
$\OO_{\PP^{n-1}}(d')$.

Now we estimate the possibilities of the integers $d'$ and $r$.
Since $Z$ is Fano, we get immediately $-n<d'<n$. We recall now that
$0<d<n$. Since $r=d'-d$, we have finally $d'<r<n+d'$. So we get the example (2)-(b). \finpreuve

\begin{prop}\label{rho}
Let $X$ be a smooth projective variety and 
let $\varphi:X\to Z$ be an elementary extremal contraction 
of ray $R$. We assume that there exists a prime divisor $E$ 
such that $\varphi(E)$ is a point. If there exists a curve 
$C\subset E$ such that $E\cdot C>0$ then $\rho(X)=1$. 
\end{prop}

{\bf Proof.} By assumption, $E\cdot R>0$ (ie: 
for all $[\Gamma]\in R$, we have $E\cdot \Gamma>0$).
If $\varphi$ is birational, $E=\Exc(\varphi)$ and 
$E\cdot R<0$
\footnote{Let $H$ be a hyperplane section of $Z$ passing through the point 
$\varphi(E)$. We can write: $\varphi^{*}H=\tilde{H}+kE$ $(k>0)$ 
where $\tilde{H}$ is the strict transform of $H$. Note that there exists 
a curve $A\subset E$ such that $\tilde{H}\cdot A>0$. We get 
$0=(\varphi^{*}H)\cdot A=\tilde{H}\cdot A+k(E\cdot A)$. Therefore 
$E\cdot A<0$. It follows that $E\cdot R<0$ because $\varphi$ is an elementary 
contraction.}, a contradiction.  
Hence $\varphi$ is of fiber type. If $\dim Z>0$, 
we take a point $z\neq\varphi(E)$. Then for any curve $B$ 
contained in $\varphi^{-1}(z)$, we get $E\cdot B=0$
although $[B]\in R$. This contradicts our assumption. 
Therefore $Z$ is a point, namely $\rho(X)=1$. \finpreuve

\section{Comments on the bound of Picard number}
By the classification result of Theorem \ref{divneg}, 
the Picard number of a Fano manifold $X$ 
containing a divisor $E\simeq \PP^{n-1}$ with $\deg N_{E/X}<0$ 
is less than or equal to $3$. 
This is in fact true in more general situation.
\begin{prop}\label{rho3}
Let $X$ be an $n$-dimensional Fano manifold $(n\geq 3)$. 
We assume that $X$ contains a prime divisor $E$ with $\rho(E)=1$. 
Then $\rho(X)\leq 3$.
\end{prop}

%This proposition has an immediate conclusion:
%\begin{cor}
%Let $X$ be a Fano manifold of dimension $n\geq 3$ and 
%$E$ be a variety with $\rho(E)=1$. Then there does not 
%exsit any morphism from $E$ to $X$ if $\rho(X)\geq 4$.
%\end{cor} 

{\bf Proof.} 
Since $X$ is Fano, we can take an extremal ray 
$\RR[f]$ such that $E\cdot f>0$. 
Let $\varphi:=\cont_{\RR[f]}: X\to Z$ be the associated 
elementary extremal contraction. 

\

 If there exists $z\in Z$ such that 
$\dim \varphi^{-1}(z)\geq 2$, then there exists a curve 
$B\subset E\cap \varphi^{-1}(z)$. We have $[B]\in\RR[f]$, 
because $B\subset\varphi^{-1}(z)$. This implies that 
$E\cdot B>0$. Since $B\subset E$, we have $\rho(X)=1$ 
by Proposition \ref{rho}.

\

 If $\dim \varphi^{-1}(z)\leq 1$ for all $z\in Z$, 
by \cite{Ando} $\varphi$ is either a conic bundle or a smooth 
blow-up of a smooth center of codimension 2 
(here, $\PP^{1}$-bundle is considered as a special case of conic bundles). 
\begin{itemize}
\item If $\varphi$ is a conic bundle (or a $\PP^{1}$-bundle), 
the restriction map $\varphi_{|E}:E\to Z$ is surjective 
(and moreover finite). 
Since $\rho(E)=1$, $\rho(Z)=1$. It follows that 
$\rho(X)=\rho(Z)+1=2$ (because $\varphi$ is elemental). 
\item We treat now the case where 
$\varphi$ is a blow-up along a 
smooth center $W$ of codimension 2. Let $E'=\varphi(E)$. 
Since there exists a surjective map 
$\varphi_{|E}:E\to E'$, we have $\rho(E')=1$.

\

{\sc Claim.} {\it The smooth variety $Z$ is Fano.}
 
{\it Proof:} By Proposition \ref{Yfano} it is sufficient to show that
for any curve $B\subset W$ we have $-K_{Z}\cdot B>0$. 
Let $A$ be a curve in $E'$ not contained in $W$. 
Since $\rho(E')=1$, there exists a positive real number $a$ such that 
$B\equiv aA$ in $E'$. Since $E'\subset Z$, this numerical equivalence 
holds also in $Z$. Therefore $-K_{Z}\cdot B>0$ 
if and only if $-K_{Z}\cdot A>0$. This is in fact the case, because 
$$
-K_{Z}\cdot A=-K_{X}\cdot\tilde{A}+F\cdot\tilde{A}>0
$$ 
where  
$\tilde{A}$ is the strict transform of $A$ by $\varphi$. 

\

So we can take an extremal ray $\RR[m]$ 
such that $E'\cdot m>0$. 
Let $\mu:=\cont_{\RR[m]}:Z\to V$  
be the  associated elementary extremal contraction. 
\begin{itemize}
\item If there exists $v\in V$ such that 
$\dim \varphi^{-1}(v)\geq 2$, by the same argument as above, 
we get $\rho(Z)=1$. Hence $\rho(X)=\rho(Z)+1=2$.
\item If $\dim \varphi^{-1}(v)\leq 1$ for all $v\in V$, 
as mentioned above, $\mu$ is either a conic bundle or 
a blow-up of codimension 2. In the first case, we have 
a surjective map $\mu_{|E'}:E'\to V$, so $\rho(V)=1$ and 
$\rho(Z)=\rho(V)+1=2$. It follows that $\rho(X)=\rho(Z)+1=3$. 

Now we rule out the birational case: let $M$ be the exceptional 
divisor of $\mu$. Since $\rho(E')=1$, $M\cap W\neq\emptyset$. 
It follows that there is a fiber $m$ of the exceptional 
divisor $M$ such that $m\cap W\neq\emptyset$.  
If $\tilde{m}$ is the strict transform of $m$ by $\varphi$, we get 
$$
K_{X}\cdot \tilde{m}=K_{Z}\cdot m+F\cdot \tilde{m}\geq -1+1=0,
$$  
a contradiction because $X$ is Fano. 

\

We conclude that in every case, we have $\rho(X)\leq 3$. 
\finpreuve

\end{itemize}     
\end{itemize}

\

\

\flushright{email: tsukiokatoru@yahoo.co.jp}
\end{document}